\documentclass[leqno,12pt]{amsart} 
\usepackage{amssymb} 
\usepackage{braket}
\usepackage{empheq}
\usepackage{bm}

\theoremstyle{plain} 
\newtheorem{theorem}{\noindent\bf Theorem}[section] 
\newtheorem{lemma}[theorem]{\noindent\bf Lemma}
\newtheorem{corollary}[theorem]{\noindent\bf Corollary}
\newtheorem{proposition}[theorem]{\noindent\bf Proposition}

\theoremstyle{definition}
\newtheorem{definition}[theorem]{\noindent\bf Definition}
\newtheorem{remark}[theorem]{\noindent\bf Remark}

\newcommand{\Ker}[0]{\operatorname{Ker}}

\newcommand{\deldel}{\sqrt{-1}\partial \overline{\partial}}
\newcommand{\dbar}{\overline{\partial}}
\newcommand{\e}{\varepsilon}

\newcommand{\lla}[0]{{\langle\!\hspace{0.02cm} \!\langle}}
\newcommand{\rra}[0]{{\rangle\!\hspace{0.02cm}\!\rangle}}

\begin{document}

\title[Calabi type functionals for coupled K\"ahler-Einstein metrics]{Calabi type functionals for coupled K\"ahler-Einstein metrics} 
\author[S. Nakamura]{Satoshi Nakamura} 


\thanks{ 
}

\address{
Department of Mathematics, 
Tokyo Institute of Technology,
2-12-1, Ookayama, Meguro-ku, Tokyo, 152-8551, Japan.
}
\email{s.nakamura@math.titech.ac.jp}

\maketitle

\begin{abstract}
We introduce the coupled Ricci-Calabi functional and the coupled H-functional 
which measure how far a K\"ahler metric is from 
a coupled K\"ahler-Einstein metric in the sense of Hultgren-Witt Nystr\"om.
We first give corresponding moment weight type inequalities
which estimate each functional in terms of algebraic invariants.
Secondly, we give corresponding Hessian formulas for these functionals at each critical point,
which have an application to a Matsushima type obstruction theorem for the existence
of a coupled K\"ahler-Einstein metric. 
\end{abstract}

 \setcounter{tocdepth}{1}
\tableofcontents

\section{Introduction}
Hultgren-Witt Nystr\"om \cite{HN18} introduced the notion of a coupled K\"ahler-Einstein metric
on a compact complex manifold of general type or a Fano manifold.
In this paper we mainly focus on the Fano case.
Let $X$ be an $n$-dimensional Fano manifold.
{\it A decomposition} of the first Chern class $2\pi c_1(X)$ is a sum
\[2\pi c_1(X)=\alpha_1 + \cdots + \alpha_N\]
where each $\alpha_i$ is a K\"ahler class for $X$.
Let $\theta_i\in\alpha_i$ be a reference K\"ahler metric and let $V_i$ be the volume $\int_X\theta_i^n$ of $\alpha_i$.
We define the set of tuples of K\"ahler potentials by
\begin{eqnarray*}
\mathcal{M}
:=
\prod_{i=1}^{N}\mathcal{M}_{i}
:=\prod_{i=1}^{N}\Set{\phi_i\in C^{\infty}(X;\mathbb{R}) | \omega_{\phi_i}
:=\theta_i+\deldel\phi_i>0},
\end{eqnarray*}
and identify a K\"ahler metric $\omega_{\phi_i}$ with its potential $\phi_i$.
A tangent space of $\mathcal{M}$ is identified with $(C^{\infty}(X;\mathbb{R}))^N$.
For any tuple of K\"ahler metrics $\Phi=(\phi_i)_{i=1}^N\in\mathcal{M}$, 
since the Ricci form $\mathrm{Ric}(\omega_{\phi_i}):=-\deldel\log\omega_{\phi_i}^n$
and the sum $\sum_{i=1}^N\omega_{\phi_i}$ represent $2\pi c_1(X)$,
there exists a unique smooth real function $f_i(\Phi)$ satisfying
\begin{equation}\label{Ricci}
{\rm Ric}(\omega_{\phi_i})-\sum_{j=1}^N\omega_{\phi_j} = \deldel f_i(\Phi) \quad\text{and}\quad \int_X(1-e^{f_i(\Phi)})\omega_{\phi_i}^n=0.
\end{equation}
In this paper we call the tuple $(f_i(\Phi))_{i=1}^N$ {\it the Ricci potential} for $\Phi=(\phi_i)_{i=1}^N$.
Then the tuple $\Phi=(\phi_i)_{i=1}^N$ is 
{\it a coupled K\"ahler-Einstein metric} for the decomposition $(\alpha_i)_{i=1}^N$
if every $f_i(\Phi)$ vanishes, that is,
$$\mathrm{Ric}(\omega_{\phi_1})=\cdots = \mathrm{Ric}(\omega_{\phi_N})=\sum_{i=1}^{N}\omega_{\phi_i}.$$

Coupled K\"ahler-Einstein metrics were studied extensively in recent years 
\cite{DP19, DH18, Fu22, FZ18, FZ19, Ha21, HN18, Hul17, Na20, Pin18, Tak19, Tak19-2}.
One of the motivation to study comes from algebro-geometric stabilities.
Indeed, 
Hultgren-Witt Nystr\"om \cite{HN18} introduced the notion of 
called K-polystability for $(X, (\alpha_i)_{i=1}^N)$, 
and showed that the existence of a coupled K\"ahler-Einstein metric implies it.
Datar-Pingali \cite{DP19} introduced a framework of geometric invariant theory 
for a coupled constant scalar curvature K\"ahler metric 
which is a generalization of a coupled K\"ahler-Einstein metric.

The well-known Calabi functional \cite{Ca82, Ca85}, which is the $L^2$-norm of a scalar curvature, 
plays an important role for studies of a K\"ahler-Einstein metric and a constant scalar curvature K\"ahler metric. 
In this paper, we introduce two Calabi type functionals which measure how $\Phi$ is far from a coupled K\"ahler-Einstein metric.
We first discuss 
moment weight type inequalities which give lower bounds of these functionals in terms of algebro-geometric stability invariants.
Secondly, we discuss Hessians for these functionals at each critical point to obtain various corollaries including a new proof of a Matsushima type obstruction theorem for the existence of a coupled K\"ahler-Einstein metric.

Let us introduce two Calabi type functionals as follows.
\[R_c(\Phi)=\sum_{i=1}^N \frac{1}{V_{i}}\int_X (1-e^{f_i(\Phi)})^2\omega_{\phi_i}^n
\quad\text{and}\quad
H_c(\Phi)=\sum_{i=1}^N \frac{1}{V_{i}}\int_X f_i(\Phi)e^{f_i(\Phi)}\omega_{\phi_i}^n.\]
In this paper we call $R_c$ {\it the coupled Ricci-Calabi functional} and
$H_{c}$ {\it the coupled H-functional}.
These are non-negative functionals in $\mathcal{M}$
whose zeros are coupled K\"ahler-Einstein metrics (see the inequality \eqref{Pinsker}).
When $N=1$, these functional are written as $R$ and $H$ respectively, and are called the Ricci-Calabi functional and the H-functional respectively. 
Functionals $R$ and $H$ were studied in 
\cite{Be16, DS16, Fo16, He16, Hi19, Na19, TZZZ13, Xi19, Ya}, 
and in particular play important roles in the context of optimal destabilizers for a Fano manifold admitting no K\"ahler-Einstein metric.

\subsection{Moment weight type inequalites}
The Calabi functional for a polarized manifold
has a lower bound in terms of the Donaldson-Futaki invariant \cite{Do05}. 
Such an inequality is called the moment weight inequality 
since it already appears in geometric invariant theory as an inequality 
between the squared norm of a moment map and a Hilbert-Mumford weight.
The Ricci-Calabi functional $R$ and the $H$-Functional $H$ satisfy
a corresponding moment weight inequality \cite{Be16, DS16, He16, Hi16, Hi19}. 
The first results in this paper are two moment weight type inequalities for $R_c$ and $H_c$
which generalize these inequalities.
\begin{theorem}\label{moment weight}
We have
$$\inf_{\Phi\in\mathcal{M}}R_c(\Phi)^{1/2} \geq 
\sup_{(\mathcal{X},(\mathcal{L}_i)_{i=1}^N)}\frac{-D_c(\mathcal{X},(\mathcal{L}_i)_{i=1}^N)}{\| (\mathcal{X},(\mathcal{L}_i)_{i=1}^N) \|_2}
\quad\text{and}\quad
\inf_{\Phi\in\mathcal{M}}H_c(\Phi) \geq
\sup_{(\mathcal{X},(\mathcal{L}_i)_{i=1}^N)}-H_c(\mathcal{X},(\mathcal{L}_i)_{i=1}^N).$$
\end{theorem}
Here $(\mathcal{X},(\mathcal{L}_i)_{i=1}^N)$ in the above supremums runs through 
arbitrary test configuration of the decomposition $(X, (\alpha_{i})_{i=1}^{N})$ introduced in \cite{HN18},
$\| (\mathcal{X},(\mathcal{L}_i)_{i=1}^N) \|_2$ is the $L^2$-norm,
$D_c(\mathcal{X},(\mathcal{L}_i)_{i=1}^N)$ is the coupled Ding-invariant,
and $H_c(\mathcal{X},(\mathcal{L}_i)_{i=1}^N)$ is the coupled $H$-invariant.
These notions are introduced in Section \ref{MWineq}.

As a direct consequence of Theorem \ref{moment weight}, 
a Fano manifold admitting a coupled K\"ahler-Einstein metric satisfies
algebraic (semi-)stability conditions.

\begin{corollary}
Suppose the existence of a coupled K\"ahler-Einstein metric
for the decomposition $(\alpha_{i})_{i=1}^{N}$ of $2\pi c_{1}(X)$.
Then we have
$$
D_{c}(\mathcal{X},(\mathcal{L}_{i})_{i=1}^{N})\geq 0
\quad\text{and}\quad
H_{c}(\mathcal{X},(\mathcal{L}_{i})_{i=1}^{N})\geq 0
$$
for any test configuration $(\mathcal{X},(\mathcal{L}_{i})_{i=1}^{N})$
for $(X,(\alpha_{i})_{i=1}^{N})$.
\end{corollary}

When $N=1$, the equalities in Theorem \ref{moment weight} in fact hold. 
Dervan-Sz\'ekelyhidi \cite{DS16} showed the moment weight equality for $H$ 
by applying the K\"ahler-Ricci flow together with deep results in \cite{CSW15, CW14}.
Hisamoto \cite{Hi19} showed corresponding equalities for $R$ and $H$ 
by using the inverse Monge-Amp\`ere flow \cite{CHT} and the K\"ahler-Ricci flow respectively, 
together with a technique for multiplier ideal sheaves.
When $N>1$, in order to establish the equality in Theorem \ref{moment weight},
 it is natural to consider the generalization of these flow, that is, the coupled inverse Monge-Amp\`ere flow
\begin{equation}\label{cKRF}
\frac{d}{dt}\phi_i(t) = 1-e^{f_i(\Phi(t))}
\quad
(i=1,\dots,N)
\end{equation}
and the coupled K\"ahler-Ricci flow
\begin{equation}\label{cIMAF}
\frac{d}{dt}\phi_i(t) = -f_i(\Phi(t))
\quad
(i=1,\dots,N).
\end{equation}
However little is known for these flows at present.
For instance the short time existence for each flow is true since they are parabolic.
However the long time existence is not established.
In Section \ref{Hessian}, we see that each flow is a gradient flow for $R_c$ and $H_c$ respectively
(Corollary \ref{flow convergense}).
They will present not only some applications to establish the equalities in Theorem \ref{moment weight}
but also some interesting problems in geometric analysis.   

\begin{remark}
Very recently, Hashimoto \cite{Ha21} introduced a different framework of test configurations
for a decomposition $(\alpha_{i})_{i=1}^{N}$
where $\alpha_{i}=2\pi c_{1}(L_{i})$ for a line bundle $L_{i}\to X$.
The author expects that corresponding moment weight type inequalities hold in his framework.
\end{remark}

\subsection{Hessian formulas for functionals and its application to a Matsushima type obstruction theorem}
In this paper, we call a critical point of $R_c$ 
{\it a coupled Mabuchi soliton} (cf. \cite{Hi19-2, Ya}) 
and call a critical point of $H_{c}$ {\it a coupled K\"ahler-Ricci soliton} (cf. \cite{He16}).
In Section \ref{Hessian} we show that 
a tuple $\Phi=(\phi_i)_{i=1}^N\in\mathcal{M}$ is a coupled Mabuchi soliton 
if and only if the vector fields $\mathrm{grad}_{\phi_i}e^{f_i(\Phi)}$ are holomorphic and
$\mathrm{grad}_{\phi_1}e^{f_1(\Phi)}=\cdots=\mathrm{grad}_{\phi_N}e^{f_N(\Phi)}$.
Similarly, we show that $\Phi=(\phi_i)_{i=1}^N\in\mathcal{M}$ is a coupled K\"ahler-Ricci soliton 
if and only if the vector fields $\mathrm{grad}_{\phi_i}f_i(\Phi)$ are holomorphic and
$\mathrm{grad}_{\phi_1}f_1(\Phi)=\cdots=\mathrm{grad}_{\phi_N}f_N(\Phi)$.

Examples of coupled Mabuchi solitons and coupled K\"ahler-Ricci solitons on Fano manifolds with large symmetry have already appeared in \cite{DH18} (see also \cite{Hul17}).
However, in that paper, the conditions 
$\mathrm{grad}_{\phi_1}e^{f_1(\Phi)}=\cdots=\mathrm{grad}_{\phi_N}e^{f_N(\Phi)}$ 
and $\mathrm{grad}_{\phi_1}f_1(\Phi)=\cdots=\mathrm{grad}_{\phi_N}f_N(\Phi)$
are not required for each definition.
Their motivation is the construction of such metrics by proving $C^{0}$-estimate 
to a class of coupled Monge-Amp\`ere equation,
which is independent of the Calabi type functionals.

The second result in this paper shows that each critical metric is in fact a local minimum of the corresponding functional 
by giving the Hessian formulas at each critical point.
Let $\lla \cdot, \cdot \rra_{\Phi}$ and $\lla\cdot, \cdot\rra^f_{\Phi}$ 
be Hermitian inner products on $(C^{\infty}(X;\mathbb{C}))^N$ defined by
\begin{equation*}
\lla \bm{u}, \bm{v} \rra_{\Phi} 
= \sum_{i=1}^N \frac{1}{V_{i}} \int_X u_i \overline{v_i}\omega_{\phi_i}^n
\quad\text{and}\quad
\lla \bm{u}, \bm{v} \rra^f_{\Phi} 
= \sum_{i=1}^N \frac{1}{V_{i}} \int_X u_i \overline{v_i}e^{f_{i}(\Phi)}\omega_{\phi_{i}}^{n}.
\end{equation*}
Let $P_{\Phi}$ and $P_{\Phi}^f$ be the operators acting on $T_{\Phi}\mathcal{M}$
defined by \eqref{Operator} and \eqref{Operatorf} in Section \ref{Hessformula}.

\begin{theorem}\label{Hessian formula}
At a coupled Mabuchi soliton $\Phi_R \in\mathcal{M}$, 
the Hessian of the coupled Ricci-Calabi functional $R_{c}$ is written as
\begin{equation}
\mathrm{Hess}(R_c)(\delta\Phi_1, \delta\Phi_2)
= 2\lla P_{\Phi_R}^f \overline{P_{\Phi_R}^f}\delta\Phi_1, \delta\Phi_2 \rra_{\Phi_R}
= 2\lla \overline{P_{\Phi_R}^f} P_{\Phi_R}^f\delta\Phi_1, \delta\Phi_2 \rra_{\Phi_R}
\end{equation}
for any variations $\delta\Phi_1, \delta\Phi_2\in T_{\Phi_R}\mathcal{M}$.
At a coupled K\"ahler-Ricci soliton $\Phi_H \in\mathcal{M}$,
the Hessian of the coupled $H$-functional is written as 
\begin{equation}
\mathrm{Hess}(H_c)(\delta\Phi_1, \delta\Phi_2)
= \lla P_{\Phi_R} \overline{P_{\Phi_R}}\delta\Phi_1, \delta\Phi_2 \rra_{\Phi_R}^f
= \lla \overline{P_{\Phi_R}} P_{\Phi_R}\delta\Phi_1, \delta\Phi_2 \rra_{\Phi_R}^f
\end{equation}
for any variations $\delta\Phi_1, \delta\Phi_2\in T_{\Phi_H}\mathcal{M}$.
\end{theorem}
Here the operator $\overline{P_{\Phi_R}} $ (resp. $\overline{P_{\Phi_R}^f}$) in the above theorem 
is the complex conjugate of $P_{\Phi_R}$ (resp. $P_{\Phi_R}^f$).
Since the operator $P_{\Phi}^f$ (resp. $P_{\Phi}$) is self-adjoint and non-negative with respective to 
$\lla \cdot, \cdot \rra_{\Phi}$ (resp. $\lla \cdot, \cdot \rra_{\Phi}^f$) (Proposition \ref{Laplacian}), 
it turns out the following.
\begin{corollary}\label{commutativity}
Operators $P_{\Phi_{R}}^f$ and $\overline{P_{\Phi_{R}}^f}$ 
(resp. $P_{\Phi_{H}}$ and $\overline{P_{\Phi_{H}}}$) are commutative.
As the result, the composition $P_{\Phi_{R}}^f\overline{P_{\Phi_{R}}^f}$ 
(resp. $P_{\Phi_{H}}\overline{P_{\Phi_{H}}}$)
is a self-adjoint non-negative operator with respect to $\lla \cdot, \cdot \rra_{\Phi_R}$
(resp. $\lla \cdot, \cdot \rra_{\Phi_H}^f$).
In particular each Hessian of $R_c$ and $H_c$ is non-negative at each critical point.
\end{corollary}

When $N=1$, the Hessian formula of the Ricci-Calabi functional at a Mabuchi soliton is obtained by the author \cite{Na19}.
On the other hand Fong \cite{Fo16} gives the Hessian formula of the H-functional at any point by a tensor calculus.
It seems to be technically difficult to apply the Fong's tensor calculus for our case where $N>1$.
A unifying technique which generalizes the author's one in \cite{Na19} 
gives the Hessian formulas for $R_c$ and $H_c$ in Theorem \ref{Hessian formula}.

By Corollary \ref{commutativity}, operators $P_{\Phi}^f$ and $\overline{P_{\Phi}^f}$
(resp. $P_{\Phi}$ and $\overline{P_{\Phi}}$) are commutative
at a coupled Mabuchi soliton (resp. at a coupled K\"ahler-Ricci soliton).
Applying this commutativity, we show a Matsushima type obstruction theorem 
for the existence 
of a coupled Mabuchi soliton and a coupled K\"ahler-Ricci soliton.

\begin{theorem}\label{Matsushima2}
Let $X$ be a Fano manifold and $\mathfrak{h}(X)$ be the space of holomorphic vector fields.
If $X$ admits a coupled Mabuchi soliton $\Phi=(\phi_i)_{i=1}^N$,
then $\mathfrak{h}(X)$ is, as a vector space, the direct sum
$$\mathfrak{h}(X)=\sum_{\lambda \geq 0} \mathfrak{h}_{\lambda}(X)$$
where $\mathfrak{h}_{\lambda}(X)$ is the $\lambda$-eigenspace of the adjoint action of the holomorphic vector field
$\mathrm{grad}_{\phi_1}(1-e^{f_1(\Phi)}) = \cdots = \mathrm{grad}_{\phi_N}(1-e^{f_N(\Phi)})$.
If $X$ admits a coupled K\"ahler-Ricci soliton $\Phi=(\phi_i)_{i=1}^N$,
then $\mathfrak{h}(X)$ has the same decomposition as above 
where
$\mathfrak{h}_{\lambda}(X)$ is the $\lambda$-eigenspace of the adjoint action of the holomorphic vector field 
$-\mathrm{grad}_{\phi_1}f_1(\Phi) = \cdots = -\mathrm{grad}_{\phi_N}f_N(\Phi)$.
Furthermore, in both cases, $\mathfrak{h}_0(X)$ coincides with the complexification of the Lie algebra of Killing vector fields for every $\omega_{\phi_i}$.
In particular, $\mathfrak{h}_0(X)$ is reductive. 
\end{theorem}

The above theorem gives a new proof of a Matsushima type obstruction theorem for 
the existence of a coupled K\"ahler-Einstein metric 
which is already proved by Hultgren-Witt Nystr\"om \cite{HN18} and Futaki-Zhang \cite{FZ18}.

\begin{corollary}\label{Matsushima1}
Let $X$ be a Fano manifold admitting a coupled K\"ahler-Einstein metric for a decomposition of $2\pi c_1(X)$.
The holomorphic automorphism group $\mathrm{Aut}(X)$ is reductive.
\end{corollary}


\subsection*{Organization}
This paper is organized as follows.
In Section \ref{MWineq}, we introduce the notion of test configurations for a decomposition
and some algebraic invariants
to prove the moment weight type inequalities for $R_{c}$ and $H_{c}$.
Some energy functionals and its slope formulas at infinity play an important role for the proof.
In Section \ref{Hessian}, we give the Hessian formulas of $R_{c}$ and $H_{c}$
to see that each critical point is a local minimum.
As an application of the Hessian formulas, Matsushima type obstruction theorems for
the existence of coupled Mabuchi solitons
and coupled K\"ahler-Ricci solitons are proved. 

\subsection*{Acknowledgments}
The author would like to thank Tomoyuki Hisamoto for helpful discussion
about geodesic rays on the space of K\"ahler metrics.
He would like to thank the referee for numerous useful suggestions
which improved the presentation of the paper.
He
is partly supported by JSPS KAKENHI Grant JP 21K20342.


\section{Moment weight type inequalities}\label{MWineq}
\subsection{Test configurations and invariants}
Following \cite{HN18}, we define the notion of test configurations 
for a decomposition $(\alpha_{i})_{i=1}^{N}$ of $2\pi c_{1}(X)$.

\begin{definition}\label{tc}
Let $L$ be an ample line bundle over a projective manifold $Y$.
A test configuration $(\mathcal{Y},\mathcal{L})$ for $(Y,L)$ is
a normal scheme $\mathcal{Y}$,
a flat surjective morphism $\mathcal{Y}\to\mathbb{C}$
and a relatively ample line bundle $\mathcal{L}$ together with
a $\mathbb{C}^{*}$-action on $\mathcal{L}$ compatible with 
the standard $\mathbb{C}^{*}$-action on $\mathbb{C}$,
such that the fiber over $1\in\mathbb{C}$ is equal to $(Y,L)$.
\end{definition}

An $\mathbb{R}$-line bundle is understood as
a formal linear combination over $\mathbb{R}$ of line bundles.
The following definition is based on the fact that
any K\"ahler class $\alpha$ on a Fano manifold can be
written as the first Chern class of an $\mathbb{R}$-line bundle,
that is,
$\alpha=c_{1}(\sum_{j}r_{j}L_{j}):=\sum_{j}r_{j}c_{1}(L_{j})$
for an ample line bundle $L_{j}$ and a positive real number $r_{j}$.

\begin{definition}\label{Ktc}
Let $\alpha$ be a K\"ahler class on a Fano manifold $X$.
A test configuration $(\mathcal{X},\mathcal{L})$ for $(X,\alpha)$ is defined as
a test configuration in the sense of Definition \ref{tc} satisfying the following.
\begin{enumerate}
\item The scheme $\mathcal{X}$ is $\mathbb{Q}$-Gorenstein.
\item There exists an $\mathbb{R}$-line bundle $L=\sum_{j}r_{j}L_{j}$ over $X$ satisfying
$\alpha=c_{1}(L)$, where each $L_{j}$ is ample and $r_{j}>0$.
\item The line bundle $\mathcal{L}$ is written as $\sum_{j}r_{j}\mathcal{L}_{j}$ where
$\mathcal{L}_{j}$ is a line bundle over $\mathcal{X}$
such that $(\mathcal{X},\mathcal{L}_{j})$ is a test configuration for $(X,L_{j})$ 
in the sense of Definition \ref{tc}. 
\end{enumerate}
\end{definition}

\begin{definition}\label{dtc}
Let $(\alpha_{i})_{i=1}^{N}$ be a decomposition of $2\pi c_{1}(X)$.
A test configuration $(\mathcal{X}, (\mathcal{L}_{i})_{i=1}^{N})$ 
for $(X, (\alpha_{i})_{i=1}^{N})$ is defined as follows.
\begin{enumerate}
\item The scheme $\mathcal{X}$ is $\mathbb{Q}$-Gorenstein.
\item For each $i$, the $\mathbb{R}$-line bundle $\mathcal{L}_{j}$ over $\mathcal{X}$
defines a test configuration $(\mathcal{X},\mathcal{L}_{j})$ for $(X,\alpha_{i})$
in the sense of Definition \ref{Ktc}.
\item The sum $\sum_{i=1}^{N}\mathcal{L}_{i}$ defines a test configuration
$(\mathcal{X},\sum_{i=1}^{N}\mathcal{L}_{i})$ for the Fano manifold $(X,-K_{X})$
in the sense of Definition \ref{tc}.
\end{enumerate}
\end{definition}

Now we define some invariants appearing in the right hand side of
moment weight type inequalities.
Gluing each $(\mathcal{X},\mathcal{L}_{i})$ with the trivial family
we have the unique $\mathbb{C}^{*}$-equivalent family
$(\bar{\mathcal{X}},\bar{\mathcal{L}_{i}})$ over $\mathbb{P}^{1}$.
This gives a compactification of a test configuration $(\mathcal{X},(\mathcal{L})_{i=1}^{N})$
for $(X,(\alpha_{i})_{i=1}^{N})$.
 Set
$$
E(\mathcal{X},\mathcal{L}_{i})=\frac{\bar{\mathcal{L}_{i}}^{n+1}}{(n+1)V_{i}}
$$
to define {\it the coupled Ding invariant} (cf. \cite{Be16}) using the log canonical threshold
\begin{eqnarray*}
D_{c}\Bigl(\mathcal{X},(\mathcal{L}_{i})_{i=1}^{N}\Bigr)
&=&
L\Bigl(\mathcal{X},\sum_{i=1}^{N}\mathcal{L}_{i}\Bigr)
-\sum_{i=1}^{N}E(\mathcal{X},\mathcal{L}_{i}) \\
&:=&
\mathrm{lct}_{(\bar{\mathcal{X}},\mathcal{B})}(\mathcal{X}|_{t=0})-1
-\sum_{i=1}^{N}\frac{\bar{\mathcal{L}_{i}}^{n+1}}{(n+1)V_{i}}
\end{eqnarray*}
where $\mathcal{X}|_{t=0}$ is the fiber over $0\in \mathbb{P}^{1}$
and $\mathcal{B}$ is the boundary divisor uniquely determined by
the properties $\mathcal{B}\sim_{\mathbb{Q}}
-K_{\bar{\mathcal{X}}/\mathbb{P}^{1}}-\sum_{i=1}^{N}\bar{\mathcal{L}_{i}}$
and $\mathrm{supp}\mathcal{B}\subset\mathcal{X}|_{t=0}$.
We can consider the $\mathbb{C}^{*}$-action on the fiber over $t=0$
to describe $E(\mathcal{X},\mathcal{L}_{i})$ in terms of the weight
$\lambda_{i,1},\dots,\lambda_{i,N_{i,k}}$ for the action on
$H^{0}(\mathcal{X}|_{t=0}, k\mathcal{L}_{i}|_{t=0})$ for fixed positive integer $k$,
where $N_{i,k}:=\dim H^{0}(\mathcal{X}|_{t=0}, k\mathcal{L}_{i}|_{t=0})$.
It is well-known that we have
$$
E(\mathcal{X},\mathcal{L}_{i})=\lim_{k\to\infty}
\frac{\sum_{j=1}^{N_{i,k}}\lambda_{i,j}}{kN_{i,k}}.
$$
In view of this formula, the equality
$E(\mathcal{X},\mathcal{L}_{i}+c\mathcal{X}|_{t=0})
=E(\mathcal{X},\mathcal{L}_{i})+c$ holds
under the constant replacing
$\mathcal{L}_{i}\mapsto\mathcal{L}_{i}+c\mathcal{X}|_{t=0}$.
Setting $\hat{\lambda_{i}}=N_{i,k}^{-1}(\lambda_{i,1}+\cdots+\lambda_{i,N_{i,k}})$
we can define the $L^{p}$-norm
$$
\|(\mathcal{X},\mathcal{L}_{i})\|_{p}=
\lim_{k\to\infty}\Bigg[ \frac{\sum_{j=1}^{N_{i,k}}| \lambda_{i,j}-\hat{\lambda_{i}} |^{p}}
{k^{p}N_{i,k}} \Bigg]^{1/p}
$$
which is invariant under the replacing 
$\mathcal{L}_{i}\mapsto\mathcal{L}_{i}+c\mathcal{X}|_{t=0}$.
We define the $L^{p}$-norm of $(\mathcal{X},(\mathcal{L})_{i=1}^{N})$ by
$$\| (\mathcal{X},(\mathcal{L}_{i})_{i=1}^{N})\|_{p}
=\Bigg[\sum_{i=1}^{N}\|(\mathcal{X},\mathcal{L}_{i})\|_{p}^{p}\Bigg]^{1/p}.$$
Finally we define {\it the coupled H-invariant} (cf. \cite{DS16})
\begin{eqnarray*}
H_{c}\Big(\mathcal{X},(\mathcal{L}_{i})_{i=1}^{N}\Big)
&=&
L\Big(\mathcal{X},\sum_{i=1}^{N}\mathcal{L}_{i}\Big)
-\sum_{i=1}^{N}F(\mathcal{X},\mathcal{L}_{i}) \\
&:=&
L\Big(\mathcal{X},\sum_{i=1}^{N}\mathcal{L}_{i}\Big)
-\sum_{i=1}^{N}\lim_{k\to\infty}
\Big[ -\log\frac{1}{N_{i,k}}\sum_{j=1}^{N_{i,k}}e^{-\frac{\lambda_{i,j}}{k}} \Big].
\end{eqnarray*}
Note that Jensen's inequality shows
$H_{c}\Big(\mathcal{X},(\mathcal{L}_{i})_{i=1}^{N}\Big)
\geq D_{c}\Big(\mathcal{X},(\mathcal{L}_{i})_{i=1}^{N}\Big)$.

\subsection{Energy functionals and geodesic rays}
We define some energy functionals on $\mathcal{M}$.
For $\Phi=(\phi_{i})_{i=1}^{N}\in\prod_{i=1}^{N}(\mathrm{PSH}(X,\theta_{i})\cap L^{\infty})$,
{\it the Monge-Amp\`ere energy} is defined by
 $$
 E_{\theta_{i}}(\phi_{i})=\frac{1}{(n+1)V_{i}}\sum_{k=0}^{n}
\int_{X}\phi_{i}
\omega_{\phi_{i}}^{k}
\wedge\theta_{i}^{n-k}.
$$
The functional $E_{\theta_{i}}$ satisfies
$E_{\theta_{i}}(\phi_{i}+C)=E_{\theta_{i}}(\phi_{i})+C$
for any $C\in\mathbb{R}$.
For $\Phi=(\phi_{i})_{i=1}^{N}\in\prod_{i=1}^{N}(\mathrm{PSH}(X,\theta_{i})\cap L^{\infty})$, set
$$
L(\Phi)=-\log\int_{X}e^{-\sum_{i=1}^{N}\phi_{i}}\theta_{0}^{n}
$$
where $\theta_{0}$ is a K\"ahler metric satisfying
$\mathrm{Ric}(\theta_{0})=\sum_{i=1}^{N}\theta_{i}$
and $\int_{X}\theta_{0}^{n}=1$.
Hultgren-Witt Nystr\"om \cite{HN18} introduced {\it the coupled Ding functional}
$$
D_{c}(\Phi)
=L(\Phi)-\sum_{i=1}^{N}E_{\theta_{i}}(\phi_{i}).
$$
For any smooth $\Phi\in\mathcal{M}$, 
by using the equality of probability measures 
\begin{equation}\label{canonical measure}
\frac{e^{-\sum_{i}\phi_{i}}\theta_{0}^{n}}{\int_{X}e^{-\sum_{i}\phi_{i}}\theta_{0}^{n}}
=\frac{e^{f_{1}(\Phi)}\omega_{\phi_{1}}^{n}}{V_{1}}
=\cdots
=\frac{e^{f_{N}(\Phi)}\omega_{\phi_{N}}^{n}}{V_{N}},
\end{equation}
we have the first variation formula
\begin{equation}\label{dD}
\delta D_{c}(\delta\Phi)=
-\sum_{i=1}^{N}\frac{1}{V_{i}}\int_{X}\delta\phi_{i}
\Bigl( 1-e^{f_{i}(\Phi)} \Bigr)\omega_{\phi_{i}}^{n}
\end{equation}
which shows that a coupled K\"ahler-Einstein metric is a critical point of $D_{c}$.

In order to relate the invariants of test configurations and the energy functional,
we introduce the notion of geodesic rays on the space of K\"ahler metrics.

\begin{definition}
Let $\theta$ is a K\"ahler form on $X$ and 
$\Delta^{*}$ be the punctured unit disc in $\mathbb{C}$.
We identify $\theta$ with its lift to $X\times\Delta^{*}$.
Let $\phi(x,\tau)$ be an upper-semicontinuous, locally integrable, $S^{1}$-invariant function 
on $X\times\Delta^{*}$ and $\Omega(x,\tau)$ be the $(1,1)$-form on $X\times\Delta^{*}$ 
defined by $\theta+\deldel\phi$.
Then $\phi^{t}(x):=\phi(x, \tau) \in\mathrm{PSH}(X,\theta)$ with $t=-\log|\tau|^{2}$ is 
{\it a subgeodesic ray} 
if the restriction of $\Omega$ to $X\times\{\tau\}$ is semipositive for all $\tau$.
Moreover $\phi^{t}(x)$ is {\it a weak geodesic ray} 
if it is a subgeodesic ray satisfying $\Omega^{n+1}=0$ on $X\times\Delta^{*}$.
\end{definition}

The optimal $C^{1,1}$-regularity of the weak geodesic ray
is proved by \cite{CTW18} (see also \cite{PS10}).

For each locally bounded weak geodesic ray $\phi_{i}^{t}\in\mathrm{PSH}(X,\theta_{i})$,
the function $t\mapsto E_{\theta_{i}}(\phi_{i}^{t})$ is affine \cite{BBGZ}.
Since $\sum_{i=1}^{N}\phi_{i}^{t}$ is a locally bounded subgeodesic ray in 
$\mathrm{PSH}(X,\sum_{i=1}^{N}\theta_{i})$, the function
$t\mapsto L(\sum_{i=1}^{N}\phi^{t}_{i})$ is convex \cite{Bern15}.
Thus the coupled Ding functional $D_{c}$ is convex along $(\phi_{i}^{t})_{i=1}^{N}$.

In view of \cite{CT08, PS07, RWN11}, a test configuration $(\mathcal{X},\mathcal{L})$ 
for a Fano manifold $X$ with an ample line bundle $L$
in the sense of Definition \ref{tc} defines a weak geodesic ray
starting from a given K\"ahler potential.
Note that it is equivalent to the rays constructed in \cite{CT08, PS07, RWN11},
since it is known the uniqueness theorem
for the completely degenerate complex Monge-Amp\`ere equation \cite{PS10}.
It follows from an argument in \cite[page 6786]{HN18} that a test configuration
$(\mathcal{X},(\mathcal{L}_{i})_{i=1}^{N})$ for 
$(X, (\alpha)_{i=1}^{N})$ in the sense of Definition \ref{dtc}
and a collection of given K\"ahler potentials $(\phi^{0}_{i})_{i=1}^{N}\in\mathcal{M}$ induce
a collection of weak geodesic rays 
$(\phi^{t}_{i})_{i=1}^{N}\in\prod_{i=1}^{N}\mathrm{PSH}(X,\theta_{i})$ for $t\geq 0$. 

The slopes at infinity of energy functionals along $(\phi^{t}_{i})_{i=1}^{N}$
play an important role in the moment weight type inequalities we are interested in.
Results of Berman \cite{Be16} showed
$$
E(\mathcal{X},\mathcal{L}_{i})
=\lim_{t\to\infty}\frac{E_{\theta_{i}}(\phi_{i}^{t})}{t}
\quad\text{and}\quad
D_{c}(\mathcal{X},(\mathcal{L}_{i})_{i=1}^{N})
=\lim_{t\to\infty}\frac{D_{c}((\phi_{i}^{t})_{i=1}^{N})}{t}.
$$ 
By the $C^{1,1}$-regularity of the weak geodesic ray $\phi_{i}^{t}$,
the existence of the time derivative $\dot{\phi^{t}_{i}}$ is guaranteed.
Berndtsson \cite{Bern18} showed that 
the push forward probability measure
$\mathrm{DH}(\mathcal{X},\mathcal{L}_{i}):=(\dot{\phi_{i}^{t}})_{*}(V_{i}^{-1}\omega_{\phi_{i}^{t}}^{n})$
on $\mathbb{R}$ is independent of $t$.
Hisamoto \cite{Hi16} showed that the weak convergence of the spectral measure
$$
\frac{1}{N_{i,k}}\sum_{j=1}^{N_{i,k}}\delta_{\frac{\lambda_{i,j}}{k}}
\to\mathrm{DH}(\mathcal{X},\mathcal{L}_{i})
$$
as $k\to\infty$
to obtain the equality
$$
\| (\mathcal{X},\mathcal{L}_{i}) \|_{p}
=\Bigg[\frac{1}{V_{i}}\int_{X}\Big| \dot{\phi_{i}^{t}}-E(\mathcal{X},\mathcal{L}_{i}) \Big|^{p}
\omega_{\phi_{i}^{t}}^{n}\Bigg]^{1/p}.
$$
Consider the ``virtual slope''
$$
F(\dot{\phi_{i}^{t}})=-\log\frac{1}{V_{i}}\int_{X}e^{-\dot{\phi_{i}^{t}}}\omega_{\phi_{i}}^{n}
=-\log\int_{\mathbb{R}}e^{-x}\mathrm{DH}(\mathcal{X},\mathcal{L}_{i})
$$
to obtain the slope formula
$$
H_{c}\Big( \mathcal{X},(\mathcal{L}_{i})_{i=1}^{N} \Big)
=\lim_{t\to\infty}\Bigg[ \frac{L(\sum_{i=1}^{N}\phi_{i}^{t})}{t}-\sum_{i=1}^{N}F(\dot{\phi_{i}^{t}}) \Bigg].
$$

\subsection{Proof of Theorem \ref{moment weight}}
\begin{proof}[Proof of the moment weight type inequality for $R_{c}$]
Fix any $\Phi=(\phi_{i})_{i=1}^{N}\in\mathcal{M}$ and
any test configuration $(\mathcal{X},(\mathcal{L}_{i})_{i=1}^{N})$ 
for a decomposition $(\alpha_{i})_{i=1}^{N}$ of $2\pi c_{1}(X)$.
Take weak geodesic rays $(\phi_{i}^{t})_{i=1}^{N}$ for $t\geq 0$ 
starting from $(\phi_{i})_{i=1}^{N}$ associated with $(\mathcal{X},(\mathcal{L}_{i})_{i=1}^{N})$.

By the convexity of the coupled Ding functional,
$$
-D(\mathcal{X},(\mathcal{L}_{i})_{i=1}^{N})
=\lim_{t\to\infty}\frac{-D_{c}((\phi_{i}^{t})_{i=1}^{N})}{t}
\leq -\frac{d}{dt}\bigg|_{t=0}D_{c}((\phi_{i}^{t})_{i=1}^{N}).
$$
By the equality \eqref{dD}, the normalization of the Ricci potentials and the Schwartz inequality,
\begin{eqnarray*}
-\frac{d}{dt}\bigg|_{t=0}D_{c}((\phi_{i}^{t})_{i=1}^{N})
&=& \sum_{i=1}^{N}\frac{1}{V_{i}}\int_{X}\dot{\phi_{i}^{0}}\Big( 1-e^{f_{i}(\Phi)} \Big)
\omega_{\phi_{i}}^{n}\\
&=&  \sum_{i=1}^{N}\frac{1}{V_{i}}\int_{X}
\Big(\dot{\phi_{i}^{0}}-E(\mathcal{X},\mathcal{L}_{i})\Big)\Big( 1-e^{f_{i}(\Phi)} \Big)
\omega_{\phi_{i}}^{n} \\
&\leq& \Bigg[ \sum_{i=1}^{N}\frac{1}{V_{i}}\int_{X}
\Big|\dot{\phi_{i}^{0}}-E(\mathcal{X},\mathcal{L}_{i})\Big|^{2}
\omega_{\phi_{i}^{0}}^{n} \Bigg]^{1/2}
\Bigg[ \sum_{i=1}^{N}\frac{1}{V_{i}}\int_{X}\Big( 1-e^{f_{i}(\Phi)} \Big)^{2}
\omega_{\phi_{i}}^{n} \Bigg]^{1/2}\\
&=&\| (\mathcal{X},(\mathcal{L}_{i})_{i=1}^{N})\|_{2} R_{c}(\Phi)^{1/2}.
\end{eqnarray*}
This completes the proof.
\end{proof}

Before a proof of the moment weight type inequality for 
the coupled H-functional $H_{c}$ we give some remarks.
For two probability measures $\mu$ and $\nu$ on $X$, 
the relative entropy is defined by
$$H(\mu | \nu)=\int_{X}\log\Big( \frac{\mu}{\nu} \Big)\mu.$$
In this terminology, $H_{c}(\Phi)$ is written as
$\sum_{i=1}^{N}H(\mu_{\Phi} | \nu_{\phi_{i}})$
where $\mu_{\Phi}$ is one of the probability measures in the equality \eqref{canonical measure},
that is,
$$
H_{c}(\Phi)
=H\Big(\frac{e^{f_{1}(\Phi)}\omega_{\phi_{1}}^{n}}{V_{1}} \Big| \nu_{\phi_{1}}\Big)
+\cdots
+H\Big(\frac{e^{f_{N}(\Phi)}\omega_{\phi_{N}}^{n}}{V_{N}} \Big| \nu_{\phi_{N}}\Big)
$$
and where $\nu_{\phi_{i}}=\omega_{\phi_{i}}^{n}/V_{i}$. 
Note that the Csisz\'ar-Kullback-Pinsker inequality yields the inequality
\begin{equation}\label{Pinsker}
\sqrt{2H\Big(\frac{e^{f_{i}(\Phi)}\omega_{\phi_{i}}^{n}}{V_{i}}\Big | \nu_{\phi_{i}}\Big)}
\geq\frac{1}{V_{i}}\int_{X}\big|1-e^{f_{i}(\Phi)}\big|\omega_{\phi_{i}}^{n}
\end{equation}
which shows that a zero point of $H_{c}$ is a coupled K\"ahler-Einstein metric.
Note also that the relative entropy has an expression 
in terms of the Legendre duality as follows \cite{Be13}.
$$
H(\mu | \nu)=\sup_{f\in C^{0}(X;\mathbb{R})}
\Big( \int_{X}f\mu-\log\int_{X}e^{f}\nu \Big).
$$

\begin{proof}[Proof of the moment weight type inequality for $H_{c}$]
We use the same notation as in the previous proof.
By using the Legendre duality expression 
and the convexity of the function $L(\sum_{i=1}^{N}\phi_{i}^{t})$ for $t\geq 0$, we have
\begin{eqnarray*}
H_{c}(\Phi)
&\geq& \sum_{i=1}^{N}\Big(\int_{X}-\dot{\phi_{i}^{0}}\mu_{\Phi}
-\log\int_{X}e^{-\dot{\phi_{i}^{0}}}\nu_{\phi_{i}}\Big)\\
&=&
-\frac{d}{dt}\Big|_{t=0}L\Big(\sum_{i=1}^{N}\phi_{i}^{t}\Big)+\sum_{i=1}^{N}F(\dot{\phi_{i}^{0}})\\
&\geq&
-\frac{L(\sum_{i=1}^{N}\phi_{i}^{t})}{t}+\sum_{i=1}^{N}F(\dot{\phi_{i}^{t}}).
\end{eqnarray*}
Taking $t\to\infty$, we get $H_{c}(\Phi)\geq -H_{c}(\mathcal{X},(\mathcal{L}_{i})_{i=1}^{N})$.
This completes the proof.
\end{proof}

\section{Hessian formulas and its application}\label{Hessian}
\subsection{Hessian formulas}\label{Hessformula}
We fix notations to obtain Hessian formulas for the coupled Ricci-Calabi functional $R_c$
and the coupled H-functional $H_c$.
For any $\Phi=(\phi_i)_{i=1}^N\in\mathcal{M}$,
we write one of the probability measures in the equality \eqref{canonical measure} as $\mu_{\Phi}$.
Let $\Delta_{\phi_i}$ be the negative Laplacian of the metric $\phi_i$,
and let $P_{\Phi}$ and $P_{\Phi}^f$ be operators acting on 
$(C^{\infty}(X;\mathbb{C}))^N$ defined by
\begin{eqnarray}\label{Operator}
P_{\Phi}(\bm{u})=
    \left(
      -\Delta_{\phi_i}u_i- \langle \dbar u_i, \dbar f_i(\Phi) \rangle_{\phi_i} 
      -\sum_{j=1}^N u_j + \int_X\sum_{j=1}^N u_j \mu_{\Phi} 
    \right)_{i=1}^{N}
\end{eqnarray}
and
\begin{eqnarray}\label{Operatorf}
P_{\Phi}^f(\bm{u})=
    \left(
      \Bigl(-\Delta_{\phi_i}u_i- \langle \dbar u_i, \dbar f_i(\Phi) \rangle_{\phi_i} 
      -\sum_{j=1}^N u_j + \int_X\sum_{j=1}^N u_j \mu_{\Phi} \Bigr) e^{f_i(\Phi)}  
    \right)_{i=1}^{N}
\end{eqnarray}
for $\bm{u}=(u_i)_{i=1}^N \in (C^{\infty}(X;\mathbb{C}))^N$.
Their complex conjugates are defined by 
$$\overline{P_{\Phi}}(\bm{u})=\overline{P_{\Phi}(\overline{\bm{u}})} 
\quad\text{and}\quad
\overline{P_{\Phi}^f}(\bm{u})=\overline{P_{\Phi}^f(\overline{\bm{u}})}.$$
Recall the Hermitian inner products $\lla \cdot, \cdot \rra_{\Phi}$ and $\lla\cdot, \cdot\rra^f_{\Phi}$
on $(C^{\infty}(X;\mathbb{C}))^N$ are defined by
\begin{equation*}
\lla \bm{u}, \bm{v} \rra_{\Phi} = \sum_{i=1}^N \int_X u_i \overline{v_i}\frac{\omega_{\phi_i}^n}{V_i}
\quad\text{and}\quad
\lla \bm{u}, \bm{v} \rra^f_{\Phi} = \sum_{i=1}^N \int_X u_i \overline{v_i}\mu_{\Phi}.
\end{equation*}
The followings are basic properties for the operator $P_{\Phi}$ 
and the inner product $\lla \cdot, \cdot \rra_{\Phi}^f$.

\begin{proposition}\label{Laplacian}
{\rm (\cite[Proposition 2.4]{Tak19-2})}
\begin{enumerate}
\item The operator $P_{\Phi}$ is self-adjoint with respect to the inner product  $\lla \cdot, \cdot \rra_{\Phi}^f$.
\item The operator $P_{\Phi}$ is non-negative, and the Kernel $\Ker P_{\Phi}$ is equal to
$$\Set{(u_i)_{i=1}^N\in (C^{\infty}(X;\mathbb{C}))^N | \mathrm{grad}_{\phi_1}u_1=\cdots =\mathrm{grad}_{\phi_N}u_N=:V \text{ and }\text{ $V$ is holomorphic}},$$
where $\mathrm{grad}_{\phi_i}u_i$ is a type $(1,0)$ gradient vector field on $X$ defined by 
$$i_{(\mathrm{grad}_{\phi_i}u_i)}\omega_{\phi_i}=\sqrt{-1}\bar{\partial}u_i.$$
\end{enumerate}
\end{proposition}
Note that the same properties as in the above proposition holds 
for $P_{\Phi}^f$ and $\lla \cdot, \cdot \rra_{\Phi}$.

We give the first variation formula of the Ricci potential to obtain that of $R_c$ and $H_c$.
Note that the variation $\delta\Phi$ of $\Phi\in\mathcal{M}$ is in 
$T_{\Phi}\mathcal{M}$ and it is identified with an element in $(C^{\infty}(X;\mathbb{R}))^N$.
\begin{lemma}\label{del Ricci}
The first variation of the Ricci potential at $\Phi=(\phi_{i})_{i=1}^{N}\in\mathcal{M}$ is given by
\begin{equation}
\delta f_i(\delta\Phi)= -\Delta_{\phi_i}\delta\phi_i -\sum_{j=1}^N\delta\phi_j + \int_X\sum_{j=1}^N\delta\phi_j \mu_{\Phi}
\end{equation}
for any variation $\delta\Phi = (\delta\phi_i)_{i=1}^N \in T_{\Phi}\mathcal{M}$.
\end{lemma}

\begin{proof}
The derivation of the first equation in \eqref{Ricci} shows
$\delta f_i(\delta\Phi)= -\Delta_{\phi_i}\delta\phi_i -\sum_{j=1}^N\delta\phi_j+C$ for some constant $C$.
The constant $C$ is equal to $\int_X\sum_{j=1}^N\delta\phi_j \mu_{\Phi}$ 
since
\begin{eqnarray*}
0
=\int_{X}\delta(e^{f_{i}(\Phi)}\omega_{\phi_{i}})/V_{i}
=\int_X\delta f_i(\delta\Phi)\mu_{\Phi}+ \int_X\Delta_{\phi_i}\delta\phi_i\mu_{\Phi}.
\end{eqnarray*}
\end{proof}

\begin{lemma}\label{1st variation}
The first variations of $R_c$ and $H_c$ at $\Phi\in\mathcal{M}$ are given by
\begin{equation*}
\delta R_c(\delta\Phi)
= 2\lla P_{\Phi}^f(e^{f_1(\Phi)}, \dots, e^{f_N(\Phi)}), \delta\Phi \rra_{\Phi}
= 2\lla \overline{P_{\Phi}^f}(e^{f_1(\Phi)}, \dots, e^{f_N(\Phi)}), \delta\Phi \rra_{\Phi}
\end{equation*}
and
\begin{equation*}
\delta H_c(\delta\Phi)
= \lla P_{\Phi}(f_1(\Phi), \dots, f_N(\Phi)), \delta\Phi \rra_{\Phi}^f
= \lla \overline{P_{\Phi}}(f_1(\Phi), \dots, f_N(\Phi)), \delta\Phi \rra_{\Phi}^f
\end{equation*}
for any variation $\delta\Phi\in T_{\Phi}\mathcal{M}$.
\end{lemma}

\begin{remark}
In lemma \ref{1st variation},
the first variation of the coupled Ricci-Calabi functional is also expressed as 
$2\lla P_{\Phi}(e^{f_1(\Phi)}, \dots, e^{f_N(\Phi)}), \delta\Phi \rra_{\Phi}^f$.
However the expression in lemma \ref{1st variation} is technically crucial for the proof of Theorem \ref{Hessian formula}.
\end{remark}

\begin{proof}
For any variation $\delta\Phi=(\delta\phi_{i})_{i=1}^{N}$, 
direct computations together with Lemma \ref{del Ricci} show
\begin{eqnarray*}
\delta\int_X (1-e^{f_i(\Phi)})^2 \frac{\omega_{\phi_i}^n}{V_i} 
&=& \int_X 2 e^{f_i(\Phi)}\delta f_i(\Phi) \mu_{\Phi}
     + \int_X e^{f_i(\Phi)}(\Delta_{\phi_i} \delta\phi_i) \mu_{\Phi} \\
&=& -\int_X e^{f_i(\Phi)}\Delta_{\phi_i}\delta\phi_i \mu_{\Phi}
      -2\int_X e^{f_i(\Phi)}\sum_{j=1}^N\delta\phi_j \mu_{\Phi} \\
 &&+2\int_X e^{f_i(\Phi)} \mu_{\Phi}\int_X\sum_{j=1}^N\delta\phi_je^{f_i(\Phi)} \mu_{\Phi}
\end{eqnarray*}
and
\begin{eqnarray*}
\delta\int_X f_i(\Phi) e^{f_i(\Phi)} \frac{\omega_{\phi_i}^n}{V_i} 
&=& \int_X \delta f_i(\delta\Phi) \mu_{\Phi} 
     + \int_X f_i(\Phi) \delta f_i(\delta\Phi) \mu_{\Phi}
     +\int_X f_i(\Phi) (\Delta_{\phi_i} \delta\phi_i) \mu_{\Phi} \\
&=& -\int_X \Delta_{\phi_i}\delta\phi_i\mu_{\Phi}
     - \int_X f_i(\Phi) \sum_{j=1}^N \delta\phi_j \mu_{\Phi}
     + \int_X f_i(\Phi) \mu_{\Phi} \int_X \sum_{j=1}^N \delta\phi_j\mu_{\Phi}.
\end{eqnarray*}
By the integration by parts, we have
$$\int_X e^{f_i(\Phi)}\Delta_{\phi_i}\delta\phi_i \mu_{\Phi}
=2\int_X\delta\phi_i(\Delta_{\phi_i}e^{f_i(\Phi)}+\langle\dbar e^{f_i(\Phi)},\dbar f_i(\Phi) \rangle_{\phi_i})\mu_{\Phi}$$
to show that
\begin{eqnarray*}
\delta R_c(\delta\Phi)
&=& 2\sum_{i=1}^N \int_X\delta\phi_i(-\Delta_{\phi_i}e^{f_i(\Phi)}-\langle\dbar e^{f_i(\Phi)},\dbar f_i(\Phi) \rangle_{\phi_i})\mu_{\Phi} \\
 && - 2\int_X\sum_{j=1}^N e^{f_j(\Phi)} \sum_{i=1}^N \delta\phi_i \mu_{\Phi}
      + 2\int_X\sum_{j=1}^N e^{f_j(\Phi)} \mu_{\Phi} \int_X \sum_{i=1}^N \delta\phi_i\mu_{\Phi} \\
&=&  2\lla P_{\Phi}^f(e^{f_1(\Phi)}, \dots, e^{f_N(\Phi)}), \delta\Phi \rra_{\Phi}.
\end{eqnarray*}
On the other hand, by the integration by parts, we have
$$\int_X\Delta_{\phi_i}\delta\phi_i\mu_{\Phi}
=\int_X\delta\phi_i(\Delta_{\phi_i}f_i(\Phi)+\langle\dbar f_i(\Phi),\dbar f_i(\Phi) \rangle_{\phi_i})\mu_{\Phi},$$
to show that
\begin{eqnarray*}
\delta H_c(\delta\Phi) 
&=& \sum_{i=1}^N \int_X\delta\phi_i(-\Delta_{\phi_i}f_i(\Phi)-\langle\dbar f_i(\Phi),\dbar f_i(\Phi) \rangle)\mu_{\Phi} \\
 && - \int_X\sum_{j=1}^N f_j(\Phi) \sum_{i=1}^N \delta\phi_i \mu_{\Phi}
      + \int_X\sum_{j=1}^N f_j(\Phi) \mu_{\Phi} \int_X \sum_{i=1}^N \delta\phi_i\mu_{\Phi} \\
&=&  \lla P_{\Phi}(f_1(\Phi), \dots, f_N(\Phi)), \delta\Phi \rra_{\Phi}^f.      
\end{eqnarray*}
Similarly, we have
$$\int_Xe^{f_i(\Phi)}\Delta_{\phi_i}\delta\phi_i\mu_{\Phi}
=2\int_X\delta\phi_i(\Delta_{\phi_i}e^{f_i(\Phi)}
+\overline{\langle\overline{\partial} e^{f_i(\Phi)},\overline{\partial} f_i(\Phi) \rangle}_{\phi_i})\mu_{\Phi}$$
and
$$\int_X\Delta_{\phi_i}\delta\phi_i\mu_{\Phi}
=\int_X\delta\phi_i(\Delta_{\phi_i}f_i(\Phi)
+\overline{\langle\overline{\partial} f_i(\Phi),\overline{\partial} f_i(\Phi) \rangle}_{\phi_i})\mu_{\Phi}$$
which yield the equalities
$$\delta R_c(\delta\Phi)
= 2\lla \overline{P_{\Phi}^f}(e^{f_1(\Phi)}, \dots, e^{f_N(\Phi)}), \delta\Phi \rra_{\Phi}$$
and
$$\delta H_c(\delta\Phi)=\lla \overline{P_{\Phi}}(f_1(\Phi), \dots, f_N(\Phi)), \delta\Phi \rra_{\Phi}^f.$$
This completes the proof.
\end{proof}

Therefore Proposition \ref{Laplacian} and Lemma \ref{1st variation} 
show that a pair $\Phi=(\phi_i)_{i=1}^N\in\mathcal{M}$ is a coupled Mabuchi soliton 
if and only if the vector fields $\mathrm{grad}_{\phi_i}e^{f_i(\Phi)}$ are holomorphic and
$\mathrm{grad}_{\phi_1}e^{f_1(\Phi)}=\cdots=\mathrm{grad}_{\phi_N}e^{f_N(\Phi)}$.
Similarly, $\Phi=(\phi_i)_{i=1}^N\in\mathcal{M}$ is a coupled K\"ahler-Ricci soliton 
if and only if the vector fields $\mathrm{grad}_{\phi_i}f_i(\Phi)$ are holomorphic and
$\mathrm{grad}_{\phi_1}f_1(\Phi)=\cdots=\mathrm{grad}_{\phi_N}f_N(\Phi)$.


Now we prove the Hessian formulas. 
The following argument generalizes the authors's one in \cite{Na19}.

\begin{proof}[Proof of Theorem \ref{Hessian formula}]
We first compute the variation $(\delta_{\Phi} P_{\Phi}^f)(e^{f_1(\Phi)}, \dots, e^{f_N(\Phi)})$ at a coupled Mabuchi soliton $\Phi=(\phi_i)_{i=1}^N\in\mathcal{M}$ 
to obtain the Hessian of the coupled Ricci-Calabi functional
\begin{eqnarray*}
\frac{1}{2} \mathrm{Hess}(R_c)(\delta\Phi,\delta\Psi)
&=& \lla \delta_{\Phi}(P_{\Phi}^f(e^{f_1(\Phi)}, \dots, e^{f_N(\Phi)})), \delta\Psi \rra_{\Phi} \\
&=& \lla (\delta_{\Phi}P_{\Phi}^f)(e^{f_1(\Phi)}, \dots, e^{f_N(\Phi)}), \delta\Psi) \rra_{\Phi} \\
&&+ \lla P_{\Phi}^f(\delta_{\Phi}(e^{f_1(\Phi)}), e^{f_2(\Phi)}, \dots, e^{f_N(\Phi)}), \delta\Psi) \rra_{\Phi} \\
&&+ \cdots \\
&&+ \lla P_{\Phi}^f(e^{f_1(\Phi)}, \dots, e^{f_{N-1}(\Phi)}, \delta_{\Phi}( e^{f_N(\Phi)})), \delta\Psi) \rra_{\Phi}
\end{eqnarray*}
where $\delta_{\Phi}$ stands for the variation along $\delta\Phi=(\delta\phi_i)_{i=1}^N$ at $\Phi$, and $\delta\Psi$ is another variation at $\Phi$. 
Now we have the holomorphic vector field $$Z_R:=\mathrm{grad}_{\phi_1}e^{f_1(\Phi)}=\cdots=\mathrm{grad}_{\phi_N}e^{f_N(\Phi)},$$
since $\Phi=(\phi_i)_{i=1}^N$ is a coupled Mabuchi soliton.
Note that $Z_R$ is also expressed as 
$$Z_R=\mathrm{grad}_{\phi_1+t \delta\phi_1}(e^{f_1(\Phi)}+t Z_R(\delta\phi_1))=\cdots=\mathrm{grad}_{\phi_N+t \delta\phi_N}(e^{f_1(\Phi)}+t Z_R(\delta\phi_N)).$$
Indeed, the equality $i_{Z_R}(\omega_{\phi_i}+t\deldel\delta\phi_i)=\sqrt{-1}\dbar(e^{f_i(\Phi)}+t Z_R (\delta\phi_i))$ holds for each $i$ and for any small $t\in(-\e,\e)$.
Set $\Phi_t=(\phi_i+t\delta\phi_i)_{i=1}^N\in\mathcal{M}$ as a perturbation of $\Phi$.
Since $Z_{R}\in\Ker P^{f}_{\Phi_{t}}$ by Proposition \ref{Laplacian}, we then have
$$P_{\Phi_t}^f( e^{f_1(\Phi)}+t Z_R(\delta\phi_1), \dots, e^{f_N(\Phi)}+t Z_R(\delta\phi_N) )=0.$$
Thus the derivative of the above equation at $t=0$ yields the equality
\begin{eqnarray}\label{del Pf}
(\delta_{\Phi} P_{\Phi}^f)(e^{f_1(\Phi)}, \dots, e^{f_N(\Phi)})
&=& - P_{\Phi}^f (Z_R(\delta\phi_1), e^{f_2(\Phi)}, \dots, e^{f_N(\Phi)}) \\ \nonumber
&& - \cdots \\ \nonumber
&& - P_{\Phi}^f (e^{f_1(\Phi)}, \dots,e^{f_{N-1}(\Phi)}, Z_R(\delta\phi_N)).
\end{eqnarray}
Therefore, by using the equation \eqref{del Pf} ,
the formula $Z_R(\delta\phi_i)= \overline{\langle \overline{\partial} \delta\phi_i, \overline{\partial}(e^{f_i(\Phi)}) \rangle}_{\phi_i}$ 
and the formula of the derivative of $e^{f_i(\Phi)}$ in Lemma \ref{del Ricci}, we obtain
\begin{eqnarray*}
&& \delta_{\Phi}(P_{\Phi}^f(e^{f_1(\Phi)}, \dots, e^{f_N(\Phi)})) \\
&=& P_{\Phi}^f (-\overline{\langle \overline{\partial} \delta\phi_1, \overline{\partial}(e^{f_1(\Phi)}) \rangle}_{\phi_1} 
                         + \delta_{\Phi}(e^{f_1(\Phi)}), 
                           \dots, 
                         -\overline{\langle \overline{\partial} \delta\phi_N, \overline{\partial}(e^{f_N(\Phi)}) \rangle}_{\phi_N} 
                         + \delta_{\Phi}(e^{f_N(\Phi)})) \\
&=& P_{\Phi}^f\overline{P_{\Phi}^f}(\delta\Phi).
\end{eqnarray*}
Similarly $\delta_{\Phi}(\overline{P_{\Phi}^f}(e^{f_1(\Phi)}, \dots, e^{f_N(\Phi)})) = \overline{P_{\Phi}^f}P_{\Phi}^f(\delta\Phi).$
This completes the proof of the Hessian formula for the coupled Ricci-Calabi functional.

In order to prove the Hessian formula for the coupled H-functional $H_c$, 
we follow the same argument as above to obtain
\begin{eqnarray}\label{del P}
(\delta_{\Phi} P_{\Phi})(f_1(\Phi), \dots, f_N(\Phi))
&=& - P_{\Phi} (Z_H(\delta\phi_1), f_2(\Phi), \dots, f_N(\Phi)) \\ \nonumber
&& - \cdots \\ \nonumber
&& - P_{\Phi} (f_1(\Phi), \dots, f_{N-1}(\Phi), Z_H(\delta\phi_N))
\end{eqnarray}
where $\Phi=(\phi_i)_{i=1}^N\in\mathcal{M}$ is a coupled K\"ahler-Ricci soliton, 
$Z_H$ is the holomorphic vector field $\mathrm{grad}_{\phi_1}f_1(\Phi)=\cdots=\mathrm{grad}_{\phi_N}f_N(\Phi)$
and $\delta\Phi:=(\delta\phi_i)_{i=1}^N$ is a variation at $\Phi$. 
By the equation \eqref{del P},
the formula $Z_H(\delta\phi_i)
=\overline{\langle \overline{\partial} \delta\phi_i, \overline{\partial}(f_i(\Phi)) \rangle}_{\phi_i}$ 
and Lemma \ref{del Ricci} we have
\begin{eqnarray*}
\mathrm{Hess}(R_c)(\delta\Phi,\delta\Psi)
&=& \lla (\delta_{\Phi}P_{\Phi})(f_1(\Phi), \dots, f_N(\Phi)), \delta\Psi) \rra_{\Phi}^f \\
&&+ \lla P_{\Phi}(\delta_{\Phi}f_1(\Phi), f_2(\Phi), \dots, f_N(\Phi)), \delta\Psi) \rra_{\Phi}^f \\
&&+ \cdots \\
&&+ \lla P_{\Phi}(f_1(\Phi), \dots, f_{N-1}(\Phi), \delta_{\Phi}f_N(\Phi)), \delta\Psi) \rra_{\Phi}^f \\
&=& \lla P_{\Phi}\overline{P_{\Phi}}(\delta\Phi), \delta\Psi \rra_{\Phi}^f,
\end{eqnarray*}
where $\delta\Psi$ is another variation at $\Phi$.
Similarly $\mathrm{Hess}(H_c)(\delta\Phi, \delta\Psi)
=\lla \overline{P_{\Phi}}P_{\Phi}(\delta\Phi), \delta\Psi \rra_{\Phi}^f$.
This completes the proof of Theorem \ref{Hessian formula}.
\end{proof}

\begin{proof}[Proof of Corollary \ref{commutativity}]
This is a consequence of Theorem  \ref{Hessian formula} 
and the non-negativity of $P_{\Phi}$ (resp. $P_{\Phi}^f$) 
with respect to $\lla \cdot, \cdot \rra_{\Phi}^f$ (resp. $\lla \cdot, \cdot \rra_{\Phi}$)
in Proposition \ref{Laplacian}.
\end{proof}

To end this subsection we discuss the coupled flows introduced in Section 1.
Lemma \ref{1st variation} shows the following.
\begin{corollary}\label{flow convergense}
The coupled Ricci-Calabi functional $R_{c}$ is monotonically decreasing along 
a coupled inverse Monge-Amp\'ere flow
in the sense of \eqref{cIMAF}.
The coupled H-functional $H_{c}$ is monotonically decreasing along
a coupled K\"ahler-Ricci flow 
 in the sense of \eqref{cKRF}.
\end{corollary}
\begin{proof}
Set $F_{R}=F_{R}(t)=(e^{f_{i}(\Phi_{R}(t))})_{i=1}^{N}$ 
where $\Phi_{R}=\Phi_{R}(t)\in\mathcal{M}$ 
is a coupled inverse Monge-Amp\'ere flow.
Set $F_{H}=F_{H}(t)=(f_{i}(\Phi_{H}(t))_{i=1}^{N}$
where $\Phi_{H}=\Phi_{H}(t)\in\mathcal{M}$ is a coupled K\"ahler-Ricci flow.
Lemma \ref{1st variation} shows
$$
\frac{d}{dt}R_{c}(\Phi_{R})=-2\lla P^{f}_{\Phi_{R}}(F_{R}), F_{R} \rra_{\Phi_{R}}
\quad\text{and}\quad
\frac{d}{dt}H_{c}(\Phi_{H})=-\lla P_{\Phi_{H}}(F_{H}), F_{H} \rra_{\Phi_{H}}^{f}
$$
which are both non-positive by Proposition \ref{Laplacian}.
This completes the proof.
\end{proof}
Since a coupled Mabuchi soliton is a self-similar solution of
a coupled inverse Monge-Amp\'ere flow $\Phi_{R}(t)$,
Corollary \ref{commutativity} suggests that 
for a Fano manifold admitting a coupled Mabuchi soliton,
this flow $\Phi_{R}(t)$ 
starting from any metric converges to a coupled Mabuchi soliton in some sense.
The same statement for a coupled K\"ahler-Ricci flow
and a coupled K\"ahler-Ricci soliton is expected to hold.
\subsection{An application to Matsushima type obstruction theorem}

As an application of the commutativity of the operators in Corollary \ref{commutativity},
we prove Theorem \ref{Matsushima2}.

\begin{proof}[Proof of Theorem \ref{Matsushima2}]
We first fix a coupled Mabuchi soliton $\Phi=(\phi_i)_{i=1}^N$.
Since operators $P_{\Phi}^f$ and $\overline{P_{\Phi}^f}$ are commutative by Corollary \ref{commutativity}, 
$\overline{P_{\Phi}^f}\in\mathrm{End}(\Ker P_{\Phi}^f)$.
Let $E_{\lambda}$ be the $\lambda$-eigenspace of $\overline{P_{\Phi}^f}$ in $\Ker P_{\Phi}^f$.
Then we have
\begin{equation}\label{decomposition}
\Ker P_{\Phi}^f = \sum_{\lambda\geq 0}E_{\lambda}.
\end{equation} 
Note that, by Proposition \ref{Laplacian}, every $V\in\mathfrak{h}(X)$ is written as 
$V=\mathrm{grad}_{\phi_1}u_1=\cdots = \mathrm{grad}_{\phi_N}u_N$ for some $(u_i)_{i=1}^N \in \Ker P_{\Phi}^f$ which are unique up to additive constants. 
Thus, setting $$\mathfrak{h}_{\lambda}(X)
=\Set{V\in\mathfrak{h}(X) | V=\mathrm{grad}_{\phi_1}u_1=\cdots = \mathrm{grad}_{\phi_N}u_N \quad\text{and}\quad (u_i)_{i=1}^N\in E_{\lambda}}$$
and using the relation \eqref{decomposition}, we have the decomposition
\begin{equation}
\mathfrak{h}(X)=\sum_{\lambda \geq 0}\mathfrak{h}_{\lambda}(X).
\end{equation}
Here we claim that $\mathfrak{h}_{\lambda}(X)$ is the $\lambda$-eigenspace of the adjoint action of 
$Z_R:=\mathrm{grad}_{\phi_1}(1-e^{f_1(\Phi)}) = \cdots = \mathrm{grad}_{\phi_N}(1-e^{f_N(\Phi)})$.
To see this, fix an element $V=\mathrm{grad}_{\phi_1}u_1=\cdots = \mathrm{grad}_{\phi_N}u_N$ in $\mathfrak{h}_{\lambda}(X)$
where $(u_i)_{i=1}^N\in E_{\lambda}$,
and observe
\begin{eqnarray*}
&&\lambda (u_1, \dots, u_N) \\
&=& \overline{P_{\Phi}^f} (u_1, \dots, u_N) \\
&=& (\overline{P_{\Phi}^f} - P_{\Phi}^f) (u_1, \dots, u_N) \\
&=& 
(-\overline{\langle \overline{\partial} \overline{u}_{1}, \overline{\partial} e^{f_1(\Phi)} \rangle}_{\phi_1}
+ \langle \dbar u_1, \dbar e^{f_1(\Phi)} \rangle_{\phi_1}, 
\dots,
-\overline{\langle \overline{\partial} \overline{u}_N, \overline{\partial} e^{f_N(\Phi)} \rangle}_{\phi_N} 
+ \langle \dbar u_N, \dbar e^{f_N(\Phi)} \rangle_{\phi_N}) \\
&=& (\{e^{f_1(\Phi)}, u_1 \}_{\phi_1}, \dots, \{e^{f_N(\Phi)}, u_N \}_{\phi_N})
\end{eqnarray*}
where $\{\cdot,\cdot\}_{\phi_i}$ is the Poisson bracket defined by 
$\{u,v\}_{\phi_i}=(\mathrm{grad}_{\phi_i}v)u-(\mathrm{grad}_{\phi_i}u)v$.
Since the map $u\to-\mathrm{grad}_{\phi_i}u$ is a complex Lie algebra homomorphism from
$(C^{\infty}(X;\mathbb{C}), \{\cdot,\cdot\}_{\phi_i})$ to $(\Gamma(TX), [\cdot,\cdot ])$ 
where $[\cdot,\cdot ]$ is the Lie bracket defined by $[Z,W]=ZW-WZ$, 
it follows that
\begin{eqnarray*}
\lambda V &=& \mathrm{grad}_{\phi_i}\lambda u_i \\
&=& \mathrm{grad}_{\phi_i}\{e^{f_i(\Phi)}, u\}_{\phi_i} \\
&=& [-\mathrm{grad}_{\phi_i}e^{f_i(\Phi)}, \mathrm{grad}_{\phi_i}u  ] \\
&=& [ Z_R, V].
\end{eqnarray*}
This shows $\mathfrak{h}_{\lambda}(X)$ is the $\lambda$-eigenspace of the adjoint action of $Z_R$.

We next focus on $\mathfrak{h}_0(X)$.
Since $E_0=\Ker P_{\Phi}^f \cap \Ker \overline{P_{\Phi}^f}$, then both the real part $(\mathrm{Re}u_i)_{i=1}^N$ 
and the imaginary part $(\mathrm{Im}u_i)_{i=1}^N$ of $(u_i)_{i=1}^N\in E_0$ are in $E_0$ again.
It follows that
$$E_0=\Set{(u_i)_{i=1}^N \in (\sqrt{-1}C^{\infty}(X;\mathbb{R}))^N | 
\mathrm{grad}_{\phi_1}u_1=\cdots =\mathrm{grad}_{\phi_N}u_N \in \mathfrak{h}(X) }\otimes\mathbb{C}$$
and thus
$$\mathfrak{h}_0(X)=
\Set{ V+\overline{V} | (u_i)_{i=1}^N \in (\sqrt{-1}C^{\infty}(X;\mathbb{R}))^N 
                                   \quad\text{and}\quad 
                                   V= \mathrm{grad}_{\phi_1}u_1=\cdots =\mathrm{grad}_{\phi_N}u_N}
\otimes\mathbb{C}.$$
According to \cite[Lemma 2.3.8]{Fubook}, the vector filed $V+\overline{V}=\mathrm{grad}_{\phi_i}u_i+\overline{\mathrm{grad}_{\phi_i}u_i}$ as above 
is killing with respect to the K\"ahler metric $\omega_{\phi_i}$.
Therefore $\mathfrak{h}_0(X)$ is the complexification of the Lie algebra of Killing vector fields for $\omega_{\phi_i}$.
This completes the proof for the case when $\Phi=(\phi_i)_{i=1}^N$ is a coupled Mabuchi soliton.

When $\Phi=(\phi_i)_{i=1}^N$ is a coupled K\"ahler-Ricci soliton, we follow the same argument as above to obtain the decomposition
$$\mathfrak{h}(X)=\sum_{\lambda \geq 0}\mathfrak{h}_{\lambda}(X),$$
where
$\mathfrak{h}_{\lambda}(X)
=\Set{V\in\mathfrak{h}(X) | V=\mathrm{grad}_{\phi_1}u_1=\cdots = \mathrm{grad}_{\phi_N}u_N \quad\text{and}\quad (u_i)_{i=1}^N\in E_{\lambda}}$
and $E_{\lambda}$ is the $\lambda$-eigenspace of $\overline{P_{\Phi}}$ in $\Ker P_{\Phi}$.
In order to finish the proof, we only check that $\mathfrak{h}_{\lambda}(X)$ is the $\lambda$-eigenspace of the adjoint action of 
$Z_H:=-\mathrm{grad}_{\phi_1}f_1(\Phi)=\cdots =-\mathrm{grad}_{\phi_N}f_N(\Phi)$.
For any $V=\mathrm{grad}_{\phi_1}u_1=\cdots = \mathrm{grad}_{\phi_N}u_N$ in $\mathfrak{h}_{\lambda}(X)$ where $(u_i)_{i=1}^N\in E_{\lambda}$, 
we have
\begin{eqnarray*}
&&\lambda (u_1, \dots, u_N) \\
&=& \overline{P_{\Phi}} (u_1, \dots, u_N) \\
&=& (\overline{P_{\Phi}} - P_{\Phi}) (u_1, \dots, u_N) \\
&=& 
(-\overline{\langle \overline{\partial} \overline{u}_1, \overline{\partial} f_1(\Phi) \rangle}_{\phi_1} 
+ \langle \dbar u_1, \dbar f_1(\Phi) \rangle_{\phi_1}, 
\dots,
-\overline{\langle \overline{\partial} \overline{u}_N, \overline{\partial} f_N(\Phi) \rangle}_{\phi_N} 
+ \langle \dbar u_N, \dbar f_N(\Phi) \rangle_{\phi_N}) \\
&=& (\{f_1(\Phi), u_1 \}_{\phi_1}, \dots, \{f_N(\Phi), u_N \}_{\phi_N}).
\end{eqnarray*}
This yields the equality $\lambda V = [Z_H, V].$
The remaining proof is very similar to the case of coupled Mabuchi solitons.
\end{proof}


\bigskip

\end{document}